\documentclass[final,letterpaper, 10 pt, conference]{ieeeconf}
\IEEEoverridecommandlockouts
\usepackage{graphicx}
\usepackage{mathtools,bbm,amsfonts}
\usepackage{multicol}
\usepackage{xcolor}

\title{\LARGE \bf Structure preserving discretization method for 1D and 2D port-Hamiltonian systems using finite differences on staggered grids}

\author{Ignacio Diaz Alastuey$^{1}$, Yann Le Gorrec$^{1}$ and Yongxin Wu$^{1}$\\
Université Marie et Louis Pasteur, SUPMICROTECH, CNRS, institut FEMTO-ST
\\ F-25000 Besançon, France\\
        {\tt\small (ignacio.diaz, legorrec, yongxin.wu)@femto-st.fr}
\thanks{This project has received funding from the European Union's Horizon Europe research and innovative programme under the Marie Skłodowska-Curie Actions (MSCA) grant agreement No. 101073558 (ModConFlex).}%
}

\begin{document}

\maketitle
\begin{abstract}
    This paper extends previous work on finite-difference schemes over staggered grids for infinite-dimensional port-Hamiltonian systems. In the one-dimensional setting, it generalizes the discretization approach originally developed for the wave equation to a broader class of systems characterized by interconnection operators that include both differential and non-differential terms, such as the Timoshenko beam equation. The paper then introduces a discretization strategy for the two-dimensional case that requires only two grids, thereby accommodating a wider range of systems, including those whose interconnection operators contain non-differential components, such as the Mindlin plate model.
\end{abstract}

\section{Introduction}
The port-Hamiltonian (PH) formalism introduced in \cite{maschke1992} has attracted growing interest over the past twenty-five years. This modelling and control framework is founded on the description of energy exchanges within multi-physical systems. It enables the treatment of a broad class of complex systems, ranging from lumped to distributed parameter systems \cite[Chapter~6]{book_vanderSchaft2017_PHS,rashad2020}; reflecting the systems passivity and their modularity, leading to results in analysis \cite{legorrec2005}, and control \cite{rodriguez2001}.

When dealing with distributed parameter systems, whether in 1D applications like the transmission line model or beam equations \cite[Chapter~7]{jacob2012}, or in 2D or 3D applications like the one presented in \cite{liu2020}; employing a discretization that preserves the port-Hamiltonian structure is highly beneficial for both simulation and control design. With this goal in mind, various strategies have emerged since \cite{golo2004-a} where a structure-preserving implementation of the mixed finite element method was first proposed. Other strategies include finite volume \cite{kotyczka2016}, partitioned finite element \cite{cardosoribeiro2018}, pseudo-spectral \cite{moulla2011}, discrete exterior calculus \cite{seslija2011} and finite differences method \cite{trenchant2018}.

Building on these previous results, we focused on the finite difference method because of its straightforward implementation. The only preliminary calculations required to obtain the discrete system are those for the interconnection matrix, and these calculations mainly depend on the neighbouring points and the interconnection differential operator. In \cite{trenchant2018}, this strategy is only implemented for the wave equation and with differential interconnection operators that do not include non-differential terms. In this paper, we aim to extend what was presented in \cite{trenchant2018} to a larger class of PH systems. We start with the 1D case, which includes systems such as the Timoshenko beam , and extend the discussion to the 2D case, encompassing similar systems defined on rectangular domains. Additionally, we demonstrate that the interconnection matrix calculation mainly depends on the mapping of the selected grid, and the resulting system is always a PH-ODE system in explicit form.

The paper is organised as follows: Section II presents the discretization strategy for 1D port-Hamiltonian systems, while Section III extends this approach to 2D systems. We then present numerical examples in Section IV, and conclude with the final remarks in Section V.
\section{Discretization of a class of 1D port-Hamiltonian systems}
Port-Hamiltonian systems (PHS) are defined by a power-preserving structure with interconnection and damping ports. This framework allows for a modular description of finite- and infinite-dimensional systems while preserving their intrinsic passivity. The advantages of working with this class of systems range from multi-physical modelling, the ability to prove well-posedness for infinite-dimensional systems, and the design of passivity-based control strategies. 

\subsection{Class of systems \label{sec:1dSystem}}
We first consider 1D distributed parameter systems where the energy is defined as
\begin{equation}
    H(t) = \int_{\Omega}\left( \mathcal{H}_q(q,\xi)+\mathcal{H}_p(p,\xi) \right)d\xi, \label{ec:1DHamiltonian}
\end{equation}
where $q$ and $p$ are generalized coordinates and $\Omega:\{\xi\in[a,b]\}$ the spatial domain of the system. One of the most common cases involves quadratic energy functions, such as kinetic energy, elastic energy from Hooke’s law, and energy stored in linear electrical components like capacitors and inductors, etc... 
We consider PDE systems on the form \footnote{In what follows we use the notations $\dot{f}(t,\xi)=\tfrac{\partial f}{\partial t}$ or $\dot{f}(t)=\tfrac{d f}{d t}$, $\partial_q {\mathcal H}=\frac{\partial {\mathcal H}}{\partial q}$ and $\partial_p {\mathcal H}=\frac{\partial {\mathcal H}}{\partial p}$}, 
\begin{equation}
    \begin{pmatrix}
        \dot{p}\\\dot{q}
    \end{pmatrix} = \mathcal{J}\begin{pmatrix}
        \frac{\partial \mathcal{H}_p}{\partial p} \\
        \frac{\partial \mathcal{H}_q}{\partial q}
    \end{pmatrix},\label{ec:FullPHS_PDE_1D}
\end{equation}
where
\[
    \mathcal{J} = P_1\frac{\partial}{\partial \xi}(\cdot)+P_0,
\]
with
\[
\begin{aligned}
	P_1 &= \begin{pmatrix}
		\mathbf{0} & \mathcal{P}_1 \\
		\mathcal{P}_1^T & \mathbf{0}
	\end{pmatrix}, &\hspace{40pt} P_0 &= \begin{pmatrix}
	\mathbf{0} & \mathcal{P}_0 \\
	-\mathcal{P}_0^T & \mathbf{0}
	\end{pmatrix},
\end{aligned}
\]
where $\mathbf{0}$ and $\mathcal{P}_0,\mathcal{P}_1$ are respectively a zero valued matrix and full rank real matrices of appropriate dimension. For the boundary conditions, we consider the boundary efforts and flows defined as (\cite{legorrec04MTNS,legorrec2005})
\[
    \begin{pmatrix}
		f_\partial \\
		e_\partial
	\end{pmatrix} = \underbrace{\frac{1}{\sqrt{2}} \begin{pmatrix}
		P_1 & -P_1 \\
		I & I
	\end{pmatrix}
	}_{R_{ext}} \tau\begin{pmatrix}
        \frac{\partial \mathcal{H}_p}{\partial p} \\
        \frac{\partial \mathcal{H}_q}{\partial q}
    \end{pmatrix},
\]
where $\tau(\cdot)$ is the trace operator. With the boundary efforts and flows variables, one can define the boundary inputs and outputs as follows
\begin{align*}
	u_\partial(t)&=W_B \begin{pmatrix}
		f_\partial \\
		e_\partial
	\end{pmatrix},&
	y_\partial(t)&=W_C \begin{pmatrix}
		f_\partial \\
		e_\partial
	\end{pmatrix},
\end{align*}
where $W=[W_B^T\;W_C^T]^T$ is invertible and satisfies $W\Sigma W^T=\Sigma$ with
\[
\Sigma = \begin{pmatrix}
	\mathbf{0} & I\\I& \mathbf{0}
\end{pmatrix},
\]
such that the energy balance reads $\dot{H}=u_\partial^Ty_\partial$. We focus on systems in which the inputs depend only on the co-energy variables given by $\partial_q \mathcal{H}_q$ or $\partial_p \mathcal{H}_p$, but not both simultaneously at the same boundary point. The outputs, in turn, are associated with the co-energy variables corresponding to the complementary state-space variables.

\subsection{Discretization}
\begin{figure}[ht]
    \centering
    \includegraphics[width=80mm]{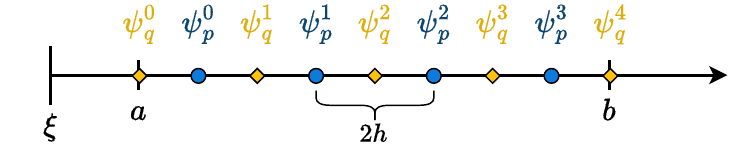}
    \caption{Example of discretization scheme.}
    \label{fig:discScheme_1D}
\end{figure}

First it is possible to notice that \eqref{ec:FullPHS_PDE_1D} can be split as
\begin{align}
    \dot{p} &= \left(\mathcal{P}_1\frac{\partial (\cdot)}{\partial \xi}+\mathcal{P}_0\right)\frac{\partial \mathcal{H}_q}{\partial q},\nonumber\\
	\dot{q} &=  \left(\mathcal{P}_1^T\frac{\partial (\cdot)}{\partial \xi}-\mathcal{P}_0^T\right)\frac{\partial \mathcal{H}_p}{\partial p}.\label{ec:dynamic_q_1D}
\end{align}
Thus the time derivative, $\dot{p}$ depends on the co-energy variable $\partial_q \mathcal{H}$, and vice-versa. With this in mind it is possible to split the approximating set of points into two different sets $\Psi_p$ and $\Psi_q$. Given that the centred differences presented in \cite{trenchant2018} are a linear approximation of the co-energy variable $\partial_p \mathcal{H}$ at points $\psi_q\in\Psi_q$ or vice-versa, the extension to non derivative term dependant on the matrix $\mathcal{P}_0$ naturally follows from a first-order Taylor polynomial approximation. As an example take the approximation of $\partial_p\mathcal{H}_p$ at a point $\psi_q\in\Psi_q$ given by two neighbouring points $\psi_p^\ell,\psi_p^{\ell+1}\in\Psi_p$. First assume that $\partial_p\mathcal{H}$ can be approximated using a first-order Taylor polynomial, and evaluate this approximation at the two given points
\[
\begin{aligned}
    \frac{\partial \mathcal{H}_p}{\partial p}(t,\psi_p^\ell) &\approx \left.\frac{\partial \mathcal{H}_p}{\partial p}\right|_{\psi_q}+\frac{\partial}{\partial \xi} \left(\frac{\partial \mathcal{H}_p}{\partial p}\right)_{\psi_q}\big(\psi_p^\ell-\psi_q\big)\\
    \frac{\partial \mathcal{H}_p}{\partial p}(t,\psi_p^{\ell+1}) &\approx \left.\frac{\partial \mathcal{H}_p}{\partial p}\right|_{\psi_q}+\frac{\partial}{\partial \xi} \left(\frac{\partial \mathcal{H}_p}{\partial p}\right)_{\psi_q}\big(\psi_p^{\ell+1}-\psi_q\big)
\end{aligned}
\]
Then considering we know the value of $\partial_p\mathcal{H}$ at each of the $\psi_p$ points we can use these two equations to solve the approximation of $\partial_p \mathcal{H}$ and its spatial derivative at $\psi_q$
\[
    \begin{aligned}
        \left.\frac{\partial \mathcal{H}_p}{\partial p}\right|_{\psi_q} \approx \frac{1}{\Delta \psi_p}\left(\Delta\psi^{\ell+1}\left.\frac{\partial \mathcal{H}_p}{\partial p}\right|_{\psi_p^{\ell}}-\Delta\psi^{\ell}\left.\frac{\partial \mathcal{H}_p}{\partial p}\right|_{\psi_p^{\ell+1}} \right) \\
        \frac{\partial}{\partial \xi}\left( \frac{\partial \mathcal{H}_p}{\partial p}\right)_{\psi_q} \approx \frac{1}{\Delta \psi_p}\left(-\left.\frac{\partial \mathcal{H}_p}{\partial p}\right|_{\psi_p^{\ell}}+\left.\frac{\partial \mathcal{H}_p}{\partial p}\right|_{\psi_p^{\ell+1}} \right)
    \end{aligned}
\]
From this, we observe that if $\psi_q$ lies at the midpoint between $\psi_p^\ell$ and $\psi_p^{\ell+1}$, then the approximation of the spatial derivative of $\partial_p \mathcal{H}$ coincides with that given by the centered difference method. Additionally, this approach ensures that the distance between the approximating points remains uniform, thus reducing the bias of the local information. With this approximation in mind it is possible to write \eqref{ec:dynamic_q_1D} at $\psi_q$ as
\[
    \dot{q}\big|_{\psi_q} \approx \frac{1}{2h}\begin{pmatrix}
    	 - \mathcal{P}_1^T -h\mathcal{P}_0^T & \mathcal{P}_1^T-h\mathcal{P}_0^T
    \end{pmatrix} \begin{pmatrix}
    	\left.\frac{\partial \mathcal{H}_p}{\partial p}\right|_{\psi_p^{\ell}} \\
    	\left.\frac{\partial \mathcal{H}_p}{\partial p}\right|_{\psi_p^{\ell+1}}
    \end{pmatrix},
\]
where $\psi_p^\ell+h=\psi_q=\psi_p^{\ell+1}-h$ leads to a centred local neighbour. The same procedure can be considered for the approximation of $\partial_q \mathcal{H}_q$ leading to a local approximation of $\dot{p}$ at a local $\psi_p$. Therefore, to generalize we can define the set of points $\Psi_q$ and $\Psi_p$ given by
\begin{align*}
	\Psi_q&:\left\{ \psi_q^m = \psi_q^0+2mh \in \Omega \big|\; m\leq M \right\}, \\
	\Psi_p&:\left\{ \psi_p^n = \psi_p^0+2nh\in \Omega \big|\; n\leq N \right\},
\end{align*}
where $\psi_q^0$ or $\psi_p^0$ is equal to $a$, $\psi_q^M$ or $\psi_p^N$ is equal to $b$, and the distance between $\psi_q^0$ and $\psi_p^0$ is $h$. An example of this partition is given by Fig. \ref{fig:discScheme_1D} where $\psi_q^0=a$, $M=4$ and $N=3$.

In the 1D case, this set of points can be organized as a vector. It is then possible to define the discretized generalized coordinate vectors $x_q=q(t,\Psi_q)$ and $x_p=p(t,\Psi_p)$, where their time derivatives are taken to be the local approximations obtained from the first-order Taylor expansion. With this, if the $m$-th element of $x_q$ is defined as $x_q^m$ its time derivative can be written as
\begin{equation}
    \dot{x}_q^m = \frac{1}{2h}\begin{pmatrix}
    	-h\mathcal{P}_0^T - \mathcal{P}_1^T & -h\mathcal{P}_0^T+\mathcal{P}_1^T
    \end{pmatrix} \begin{pmatrix}
    	\left.\frac{\partial \mathcal{H}_p}{\partial p}\right|_{\psi_p^{n_0}} \\
    	\left.\frac{\partial \mathcal{H}_p}{\partial p}\right|_{\psi_p^{n_1}}
    \end{pmatrix},\label{ec:dynamic_xq_1D}
\end{equation}
where $\psi_p^{n_0}<\psi_q^m<\psi_p^{n_1}$ and $\tfrac{\partial\mathcal{H}_p}{\partial p}\big|_{\psi_p^{n}}=\tfrac{\partial\mathcal{H}_p}{\partial p}(x_p^{n},\psi_p^{n})$. Similarly if $x_p^n$ is the $n$-th element of $x_p$ ts time derivative can be written as 
\begin{equation}
    \dot{x}_p^n = \frac{1}{2h}\begin{pmatrix}
		h\mathcal{P}_0 - \mathcal{P}_1 & h\mathcal{P}_0+\mathcal{P}_1
	\end{pmatrix} \begin{pmatrix}
	\left.\frac{\partial \mathcal{H}_q}{\partial q}\right|_{\psi_q^{m_0}} \\
	\left.\frac{\partial \mathcal{H}_q}{\partial q}\right|_{\psi_q^{m_1}}
	\end{pmatrix},\label{ec:dynamic_xp_1D}
\end{equation}
where $\psi_q^{m_0}<\psi_p^n<\psi_q^{m_1}$ and $\tfrac{\partial\mathcal{H}_q}{\partial q}\big|_{\psi_q^{m}}=\tfrac{\partial\mathcal{H}_q}{\partial q}(x_q^{m},\psi_q^{m})$.

Now we can combine \eqref{ec:dynamic_xq_1D} and \eqref{ec:dynamic_xp_1D} for each discretized generalized coordinate in $x_q$ and $x_p$. Defining the discretized state space as $x=(x_p^T\;\;x_q^T)^T$ it can be observed that the dynamic equations take the following structure:
\begin{equation}
	\dot{x} = \begin{pmatrix}
		\mathbf{0} & \frac{1}{2h}\mathcal{P}_{dp}\\ 
		\frac{1}{2h}\mathcal{P}_{dq} & \mathbf{0}
	\end{pmatrix}\begin{pmatrix}
		\frac{\partial \mathcal{H}_p}{\partial p}(x_p,\psi_p) \\
		\frac{\partial \mathcal{H}_q}{\partial q}(x_q,\psi_q)
	\end{pmatrix},\label{ec:1D_Discrete}
\end{equation}
where 
\[
\begin{aligned}
    \mathcal{P}_{dp}&=\mathcal{I}_1^{(p)}\otimes\mathcal{P}_1+\mathcal{I}_0^{(p)}\otimes\mathcal{P}_0h,\\
    \mathcal{P}_{dq}&=\mathcal{I}_1^{(q)}\otimes\mathcal{P}_1^T-\mathcal{I}_0^{(q)}\otimes\mathcal{P}_0^Th,
\end{aligned}
\]
Here $\mathcal{I}_i^{(j)}$ represents the connection coefficients associated with the $\mathcal{P}_i$ matrix for the points in the $\Psi_j$ set and $\otimes$ denotes the Kronecker product between the two matrices. To preserve a port-Hamiltonian structure, we need to define the discretized Hamiltonian and specify appropriate input/output ports. To define the discretized Hamiltonian we can use a Riemann sum structure given by
\begin{equation}
	H_d(x)=\left(\sum_{i=1}^{M^*} \mathcal{H}_q (x_q^i,\psi_q^i) 2h + \sum_{i=1}^{N^*} \mathcal{H}_p (x_p^i,\psi_q^i) 2h\right),\label{ec:1DHamiltonian_Discrete}
\end{equation}
where $M^*$ and $N^*$ are the last index in $\Psi_q$, or $\Psi_p$ respectively, that is not in the boundary of the domain. This double Riemann sum converges to the Hamiltonian defined in \eqref{ec:1DHamiltonian} if the energy density is Riemann integrable and the limit of the intervals goes to zero. From \eqref{ec:1D_Discrete} and \eqref{ec:1DHamiltonian_Discrete} we can write
\[
    \dot{x} = \begin{pmatrix}
		\mathbf{0} & \frac{1}{4h^2}\mathcal{P}_{dp}\\ 
		\frac{1}{4h^2}\mathcal{P}_{dq} & \mathbf{0}
	\end{pmatrix}\begin{pmatrix}
		\frac{\partial H_d}{\partial x_p} \\
		\frac{\partial H_d}{\partial x_q}
	\end{pmatrix},
\]

Finally, when defining the ports, we must ensure that the energy of the discrete system is conserved. Let us consider the inputs corresponding to the discrete generalized coordinates $x_q$ is the co-energy variable $\partial_p\mathcal{H}_p$. If we consider the elemental $x_q$ domain as $a=\psi_p^{n}$ and $b=\psi_p^{n+1}$, where $\psi_p^n<\psi_q^{e}<\psi_p^{n+1}$, the energy within this domain is determined solely by $x_q^e$, and the co-energy variables in $a$ and $b$ are defined as inputs. With this we take the following Hamiltonian
\[
    H_{dq}^e(t) = 2h\mathcal{H}_q(x_q^e,\psi_q^e),
\]
Then a port-Hamiltonian input-output formulation is given by
\begin{equation}
	\begin{aligned}
		\dot{x}_q^{e} &= 0 \frac{\partial H_{dq}^e}{\partial x_q^e} + B_{dq}\begin{pmatrix}
			u_p^{n} \\
			u_p^{n+1}
		\end{pmatrix},\\
		\begin{pmatrix}
			y_q^{n}\\
			y_q^{n+1}
		\end{pmatrix}&= B_{q}^T\frac{\partial H_{dq}^e}{\partial x_q^e},
	\end{aligned}
	\label{ec:elementalQ}
\end{equation}
where
\[
    B_{q}=\frac{1}{2h}\begin{pmatrix}
		-\mathcal{P}_1^T - h\mathcal{P}_0^T & \mathcal{P}_1^T - h\mathcal{P}_0^T
	\end{pmatrix}.
\]
It can be observed that, in the mechanical domain analogy, if the input corresponds to the generalized velocity, $\partial_p\mathcal{H}_p$, the associated forces represent the integral of the linearized decay of the local stress influencing the given generalized velocity over a $2h$ interval. Defining the input of an elemental $x_p$ domain, with $a=\psi_q^{m}$ and $b=\psi_q^{m+1}$, where $\psi_q^n<\psi_p^{e}<\psi_q^{n+1}$, as the power conjugated output of the elemental $x_q$ domain, the system can then be expressed using the energy in this domain as
\[
H_{dp}^e(t) = 2h\mathcal{H}_p(x_p^e),
\]
Then a port-Hamiltonian input-output formulation is given by
\begin{equation}
	\begin{aligned}
		\dot{x}_p^{e} &= 0 \frac{\partial H_{dp}^e}{\partial x_p^e} + \frac{1}{2h}\begin{pmatrix}
			u_q^{m} \\
			u_q^{m+1}
		\end{pmatrix},\\
		\begin{pmatrix}
			y_p^{m}\\
			y_p^{m+1}
		\end{pmatrix}&= \frac{1}{2h}\frac{\partial H_{dq}^e}{\partial x_q^e},
	\end{aligned}
	\label{ec:elementalP}
\end{equation}
where the conjugated output is given by $\partial_p\mathcal{H}_p$ and therefore can be used to interconnect with an $x_q$ domain.

\subsection{Energy and Structure preserving interconnection \label{ssec:EnergyStructurePreservation}}
We consider now the interconnection of two subsystems. Without loss of generality, we consider subsystem 1, defined on the spatial domain with boundaries $a_1$ and $b_1$. The input at $\xi=a_1$ is the co-energy variable associated with $q$ while the input at $\xi=b_1$ corresponds to the co-energy variable associated with $p$. Similarly, subsystem 2 is defined on the interval with boundaries $a_2$ and $b_2$, with $b_1=a_2$. The input at $\xi=a_2$ is the co-energy variable associated with $q$ and the input at $\xi=b_2$ s the co-energy variable associated with $p$. Using \eqref{ec:1D_Discrete}, \eqref{ec:1DHamiltonian_Discrete}, \eqref{ec:elementalQ} and \eqref{ec:elementalP} the subsystems can be written as follows:
\begin{equation}
    \begin{aligned}
		\begin{pmatrix}
			\dot{x}_{ip}\\
			\dot{x}_{iq}
		\end{pmatrix} &= \begin{pmatrix}
			\mathbf{0} & \mathcal{P}_{dp}^{(i)} \\
			\mathcal{P}_{dq}^{(i)} & \mathbf{0}
		\end{pmatrix}\begin{pmatrix}
		\frac{\partial H_{dp}^{(i)}}{\partial x_{ip}} \\
		\frac{\partial H_{dq}^{(i)}}{\partial x_{iq}}
		\end{pmatrix}+B_d^{(i)} \begin{pmatrix}
		u_q^{(a_i)}\\
		u_p^{(b_i)}
		\end{pmatrix} \\
		\begin{pmatrix}
			y_p^{(a_i)}\\
			y_q^{(b_i)}
		\end{pmatrix} &= {B_d^{(i)}}^T\begin{pmatrix}
		\frac{\partial H_{dp}^{(i)}}{\partial x_{ip}} \\
		\frac{\partial H_{dq}^{(i)}}{\partial x_{iq}}
		\end{pmatrix}
	\end{aligned}\label{ec:discretised_ODE_1D}
\end{equation}
where for convenience we define $x_{ip}$ and $x_{iq}$ as the vectors of discretized states $x_p^n$ and $x_q^m$ respectively ordered by their index, and where
\[
    B_d^{(i)} = \begin{pmatrix}
        B_p^{(i)} &  \mathbf{0}\\
	    \mathbf{0} & B_q^{(i)}
	\end{pmatrix},
\]
using the previously defined input matrix for elemental domains. Then, by applying a power-preserving interconnection at the interface 
$b_1=a_2$, the two subsystems can be coupled consistently:
\[
\begin{aligned}
	u_p^{b_1} &= y_p^{a_2}\\
	u_q^{a_2} &= -y_q^{b_1}\\
\end{aligned}
\]
The complete system can be expressed as 
\[
\begin{aligned}
    \dot{x}_d = \begin{pmatrix}
		\mathbf{0} & \mathcal{P}_{dp}^{(12)}\\
		  \mathcal{P}_{dq}^{(12)}  & \mathbf{0}
	\end{pmatrix} \frac{\partial H_{d}}{\partial x_d}^{(12)}+ \begin{pmatrix}
		B_p^{(1)} & \mathbf{0} \\
		\mathbf{0} & \mathbf{0} \\
		\mathbf{0} & \mathbf{0} \\
		\mathbf{0} & B_q^{(2)}
	\end{pmatrix}\begin{pmatrix}
		u_q^{a_1}\\
		u_p^{b_2}
	\end{pmatrix},\\
    \begin{pmatrix}
    	y_p^{(a_1)}\\y_q^{(b_2)} 
    \end{pmatrix} = \begin{pmatrix}
    	{B_q^{(1)}}^T & \mathbf{0} &
    	\mathbf{0} & \mathbf{0} \\
    	\mathbf{0} & \mathbf{0} &
    	\mathbf{0} & {B_p^{(2)}}^T
    \end{pmatrix}\frac{\partial H_{d}}{\partial x_d}^{(12)},
\end{aligned}
\]
where the index $(\cdot)^{(12)}$ refers to the interconnected system, with
\[
\begin{aligned}
    H_d^{(12)}(x_d) &= H_d^{(1)}(x_{1p},x_{1q})+ H_d^{(2)}(x_{2p},x_{2q})\\
    \mathcal{P}_{dp}^{(12)}&=\begin{pmatrix}
        \mathcal{P}_{dp}^{(1)} & \mathbf{0} \\
        B_q^{(2)}{B_p^{(1)}}^T & \mathcal{P}_{dp}^{(2)}
    \end{pmatrix},\\
    \mathcal{P}_{dq}^{(12)}&=\begin{pmatrix}
        \mathcal{P}_{dq}^{(1)} & -B_p^{(1)}{B_q^{(2)}}^T \\
        \mathbf{0} & \mathcal{P}_{dq}^{(2)}
    \end{pmatrix},
\end{aligned}
\]
and $x_d=\Big({x_p^{(1)}}^T\;{x_p^{(2)}}^T\;{x_q^{(1)}}^T\;{x_q^{(2)}}^T\Big)^T$. We observe that the new blocks $\mathcal{P}_j^{(12)}$ preserve the same structural form as $\mathcal{P}_j^{(i)}$. Furthermore, since the overall system can be constructed from elemental domains interconnected recursively, each interconnection matrix block is derived from the skew-symmetric blocks of the input mapping matrix. Therefore $\mathcal{P}_{dp}^{(12)}=-{\mathcal{P}_{dq}^{(12)}}^T$ and the complete system has a PH structure, with an energy variation given by the inner product between the input vector and output vectors.

\section{Extension to 2D port-Hamiltonian systems}
One of the main challenges in extending results from the 1D case to the 2D case arises from the distributed nature of the boundary. Nevertheless, the port-Hamiltonian (PH) formulation remains applicable.

\subsection{Class of systems}
We now consider a 2D distributed-parameter system whose energy is described in a manner analogous to \eqref{ec:1DHamiltonian} but with the domain defined over a two-dimensional spatial region
\begin{equation}
	H(t) = \iint_{\Omega}\Big\{ \mathcal{H}_q(q)+\mathcal{H}_p(p) \Big\}dA, \label{ec:2DHamiltonian}
\end{equation}
where $q$ and $p$ are generalized coordinates and $\Omega\subset\mathbb{R}^2$ is the spatial domain on which is defined the system. Again, one of the most common cases is that of quadratic energy functions. Therefore the formulation in \eqref{ec:FullPHS_PDE_1D} can be extended to account for the two spatial coordinates:
\begin{equation}
    \begin{pmatrix}
        \dot{p} \\ \dot{q}
    \end{pmatrix} = \mathcal{J}\begin{pmatrix}
        \frac{\partial \mathcal{H}_p}{\partial p} \\
        \frac{\partial \mathcal{H}_q}{\partial q}
    \end{pmatrix},\label{ec:FullPHS_PDE_2D}
\end{equation}
where, just like in the 1D case, the $\mathcal{J}$ operator can be summarised by
\[
    \mathcal{J} = P_1\frac{\partial}{\partial \xi_1}(\cdot)+P_2\frac{\partial}{\partial \xi_2}(\cdot)+P_0,
\]
where
\[
	P_i = \begin{pmatrix}
		\mathbf{0} & \mathcal{P}_i \\
		\mathcal{P}_i^T & \mathbf{0}
	\end{pmatrix}\;|\;i\in{1,2} \hspace{30pt} P_0 = \begin{pmatrix}
		\mathbf{0} & \mathcal{P}_0 \\
		-\mathcal{P}_0^T & \mathbf{0}
	\end{pmatrix},
\]
and $\mathbf{0}$ is a zero valued matrix of appropriate dimension. Accordingly, the time derivative of the energy is given by
\[
	\dot{H} = \oint_{\partial \Omega} \begin{pmatrix}
		\frac{\partial \mathcal{H}_p}{\partial p}^T\mathcal{P}_1\frac{\partial \mathcal{H}_q}{\partial q}^T & \frac{\partial \mathcal{H}_p}{\partial p}^T\mathcal{P}_2\frac{\partial \mathcal{H}_q}{\partial q}^T
	\end{pmatrix}\cdot \hat{n}\; ds,
\]
where $\hat{n}$ is the outward-pointing unit normal vector on the boundary. 
With this in mind, the inputs and outputs can be expressed as linear combinations of the co-energy variables defined on the system’s boundary. Again we are interested in systems where the inputs depends only on the co-energy variables of either $q$ or $p$, but not both at the same time in the same boundary point, while the output is related to the co-energy variables of the remaining state space variables.

\subsection{Discretization}
\begin{figure}[ht]
    \centering
    \includegraphics[width=80mm]{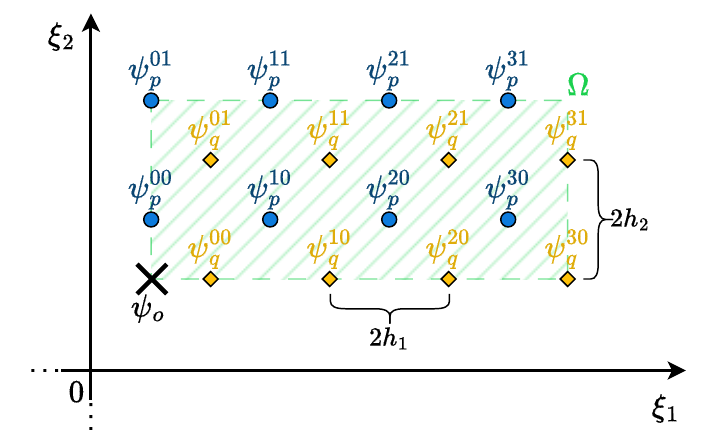}
    \caption{Example of discretization scheme.}
    \label{fig:discScheme_2D}
\end{figure}

By analogy with the 1D case, \eqref{ec:FullPHS_PDE_2D} can also be split into
\begin{align}
    \dot{p} &= \left[\mathcal{P}_1\frac{\partial (\cdot)}{\partial \xi_1}+\mathcal{P}_2\frac{\partial (\cdot)}{\partial \xi_2}+\mathcal{P}_0\right]\frac{\partial \mathcal{H}_q}{\partial q},\nonumber\\
	\dot{q} &= \left[\mathcal{P}_1^T\frac{\partial (\cdot)}{\partial \xi_1}+\mathcal{P}_2^T\frac{\partial (\cdot)}{\partial \xi_2}-\mathcal{P}_0^T\right]\frac{\partial \mathcal{H}_p}{\partial p}.\label{ec:dynamic_q_2D}
\end{align}
Thus as in the previous case it is appropriate to partition the approximating set of points into two distinct sets $\Psi_p$ and $\Psi_q$. For the remaining part of this paper, we will consider rectangular grids for the 2D case as it is the simpler case. With this in mind the Taylor's polynomial needs to include a second order derivative term. As a first approach we propose to use the crossed derivative term thus the approximation is given by
\begin{equation}
    \begin{aligned}
        f(\xi_1,\xi_2) \approx f(a,b) + (\xi_1-a)\frac{\partial f}{\partial \xi_1}(a,b)+\\
        (\xi_2-b)\frac{\partial f}{\partial \xi_2}(a,b)+(\xi_1-a)(\xi_2-b)\frac{\partial^2 f}{\partial \xi_1\partial \xi_2}(a,b).
    \end{aligned}\label{ec:Taylor2D}
\end{equation}
Following the same reasoning as in the 1D case, we define the sets of points $\Psi_q$, $\Psi_p$ which are subsets of $\Omega$ given by
\begin{align*}
	\Psi_q&:\Big\{\psi_q^{mn}=\Big(h_1\cdot(2m+m_q)\;,\;h_2\cdot(2n+n_q)\Big)+\psi_o\Big\},\\
	\Psi_p&:\Big\{\psi_p^{mn}=\Big(h_1\cdot(2m+m_p)\;,\;h_2\cdot(2n+n_p)\Big)+\psi_o\Big\},
\end{align*}
where $m_q$, $m_p$, $n_q$ and $n_p\in\{0,1\}$, with
$m_q+m_p=1$ and $n_q+n_p=1$, the pair $(m,n)$ are integers, and $\psi_o$ is an offset point. An example of this partition is given in Fig. \ref{fig:discScheme_2D} where $m_q=n_p=1$ and $m_p=n_q=0$. In this example, there are three $\psi_q$ and three $\psi_p$ points located inside $\Omega\setminus\partial\Omega$ and five of each positioned on the boundary.

We can consider the set of discrete general coordinates
\begin{align*}
	\mathcal{X}_q&:\Big\{x_q^{mn}=q(t,\psi_q^{mn})|\;\forall \psi_q^{mn}\in\Psi_q\setminus \partial \Omega\Big\} ,\\
	\mathcal{X}_p&:\Big\{x_p^{mn}=p(t,\psi_p^{mn})|\;\forall \psi_p^{mn}\in\Psi_p\setminus \partial \Omega\Big\} ,
\end{align*}
Then, by taking the time derivative of an element in $\mathcal{X}_q$ and applying the approximation from \eqref{ec:Taylor2D} to \eqref{ec:dynamic_q_2D} , in the same manner as was done for the 1D case with \eqref{ec:dynamic_xq_1D}, we take a $\psi_q^{k\ell}$ at the centre of the polygon made by the following neighbouring points $\psi_p^{mn}$, $\psi_p^{(m+1)n}$, $\psi_p^{m(n+1)}$ and $\psi_p^{(m+1)(n+1)}$ and we obtain
\begin{equation}
    \dot{x}_q^{k\ell} = \frac{1}{4h_1h_2} \mathcal{P}_{dq}^{(mn)} \frac{\partial \mathcal{H}_p}{\partial x_p^{loc}},\label{ec:dynamic_xq_2D}
\end{equation}
where
\[
    \mathcal{P}_{dq}^{(mn)} = \mathcal{I}_1\otimes\mathcal{P}_1^Th_2+\mathcal{I}_2\otimes\mathcal{P}_2^Th_1-\mathcal{I}_0\otimes\mathcal{P}_0^Th_1h_2,
\]
the $\mathcal{I}_i$ elements represent the mapping to the neighbouring points, and
\[
    x_p^{loc}=\begin{pmatrix}
        x_p^{mn}\\{x_p^{(m+1)n}}\\{x_p^{m(n+1)}}\\{x_p^{(m+1)(n+1)}}
    \end{pmatrix}.
\]
Similarly, by considering a point $\psi_p^{k\ell}$ at the center of the polygon defined by the neighboring points $\psi_q^{mn}$, $\psi_q^{(m+1)n}$, $\psi_q^{m(n+1)}$ and $\psi_q^{(m+1)(n+1)}$ , we can derive an equation analogous to \eqref{ec:dynamic_xp_1D}. If we additionally define the discretized Hamiltonian, we can employ a Riemann sum structure similar to  \eqref{ec:1DHamiltonian_Discrete}, but using an area of $4h_1h_2$ which leads to
\begin{equation}
	\begin{aligned}
	    H_d(\mathcal{X}_q,\mathcal{X}_p)=\left(\sum_{m\in\mathcal{M}_q}\sum_{n\in\mathcal{N}_q} \mathcal{H}_q (x_q^{mn},\psi_q^{mn}) +\right.\\\left.\sum_{m\in\mathcal{M}_p}\sum_{n\in\mathcal{N}_p} \mathcal{H}_p (x_p^{mn},\psi_p^{mn}) \right)4h_1h_2,
	\end{aligned}\label{ec:2DHamiltonian_Discrete}
\end{equation}
where $\mathcal{M}_q$ and $\mathcal{N}_q$ are the sets of respective indices where $\forall(m,n)\in\mathcal{M}_q\times\mathcal{N}_q\implies\psi_q^{mn}\in\Psi_q\setminus \partial \Omega$, while the same is true for $\mathcal{M}_p$ and $\mathcal{N}_p$.

To obtain an equation similar to \eqref{ec:discretised_ODE_1D}, it is important to order the discrete general coordinates into a vector. With this in mind, we incorporate each element of the sets $\mathcal{X}_p$ and $\mathcal{X}_q$ into the ordered vectors $x_p$ and $x_q$ respectively. Then by defining the discretized state space as $x=(x_p^T\;\;x_q^T)^T$ combined with \eqref{ec:2DHamiltonian_Discrete} and using an analogous input/output formulation than in the 1D case using an ordered vector for the different boundary points we obtain the dynamic equations 
\begin{equation}
    \begin{aligned}
		\begin{pmatrix}
			x_{p}\\
			x_{q}
		\end{pmatrix} &= \begin{pmatrix}
			\mathbf{0} & \mathcal{P}_{dp} \\
			\mathcal{P}_{dq} & \mathbf{0}
		\end{pmatrix}\begin{pmatrix}
		\frac{\partial H_{dp}}{\partial x_p} \\
		\frac{\partial H_{dq}}{\partial x_q}
		\end{pmatrix}+B_d \begin{pmatrix}
		u_q\\
		u_p
		\end{pmatrix} \\
		\begin{pmatrix}
			y_p\\
			y_q
		\end{pmatrix} &= {B_d}^T\begin{pmatrix}
		\frac{\partial H_{dp}}{\partial x_p} \\
		\frac{\partial H_{dq}}{\partial x_q}
		\end{pmatrix}
	\end{aligned}\label{ec:discretised_ODE_2D}
\end{equation}
where 
\[
\begin{aligned}
    \mathcal{P}_{dp}&=\mathcal{I}_1^{(p)}\otimes\mathcal{P}_1h_2+\mathcal{I}_2^{(p)}\otimes\mathcal{P}_2h_1+\mathcal{I}_0^{(p)}\otimes\mathcal{P}_0h_1h_2,\\
    \mathcal{P}_{dq}&=\mathcal{I}_1^{(q)}\otimes\mathcal{P}_1^Th_2+\mathcal{I}_2^{(q)}\otimes\mathcal{P}_2^Th_1-\mathcal{I}_0^{(q)}\otimes\mathcal{P}_0^Th_1h_2,
\end{aligned}
\]
with $\mathcal{I}_i^{(j)}$ representing the connection coefficients related to the $\mathcal{P}_i$ matrix for the points in the $\Psi_j$ set. And where
\[
    B_d = \frac{1}{4h_1h_2}\begin{pmatrix}
		B_p & \mathbf{0} \\
		\mathbf{0} & B_q 
	\end{pmatrix},
\]
with
\[
B_p = \mathcal{I}^{(up)}\otimes I_n,
\]
\[
B_q = \left(\mathcal{I}_1^{(u)}\otimes\mathcal{P}_1^Th_2+\mathcal{I}_2^{(u)}\otimes\mathcal{P}_2^Th_1-\mathcal{I}_0^{(u)}\otimes\mathcal{P}_0^Th_1h_2\right),
\]
where $\mathcal{I}^{(up)}$ represents the mapping from $\Psi_q \cap \partial \Omega$ to $\Psi_p$, $n$ is the dimension of the generalized coordinate $q$, and $\mathcal{I}_i^{(u)}$ represents the mapping from $\Psi_p \cap \partial \Omega$ to $\Psi_q$. Finally, from \eqref{ec:discretised_ODE_2D} and \ref{ssec:EnergyStructurePreservation} it is possible to extend that the discretized system has a PH structure and the discrete energy variation is given by the inner product between the input vector and the output vector.

\section{Numerical Examples}
This section presents selected numerical examples illustrating the 1D and 2D cases.
\subsection{Example: Timoshenko beam application}
For the 1D case, we consider the discretization of the Timoshenko beam equation given in PH form in \cite{ponce2024}.
We consider a beam clamped at $a=0$ and subjected to a force at the free end $b=L$. For the discretization we take $M=N=5$, which yields a step size of $h=\tfrac{L}{2\cdot5+1}$. With this, the connection coefficients in \eqref{ec:1D_Discrete} are given by
\begin{align}
    \mathcal{I}_1^{(p)} = -\left(\mathcal{I}_1^{(q)}\right)^T &= \begin{pmatrix}
        -1 & 1 & 0 & 0 & 0 \\
        0 & -1 & 1 & 0 & 0 \\
        0 & 0 & -1 & 1 & 0 \\
        0 & 0 & 0 & -1 & 1 \\
        0 & 0 & 0 & 0 & -1 
    \end{pmatrix},\label{ec:I1_p_example}\\
    \mathcal{I}_0^{(p)} = \left(\mathcal{I}_0^{(q)}\right)^T &= \begin{pmatrix}
        1 & 1 & 0 & 0 & 0 \\
        0 & 1 & 1 & 0 & 0 \\
        0 & 0 & 1 & 1 & 0 \\
        0 & 0 & 0 & 1 & 1 \\
        0 & 0 & 0 & 0 & 1 
    \end{pmatrix}.\label{ec:I0_p_example}
\end{align}
Finally, the input matrix is given by
\[
    \begin{aligned}
        B_p &= \frac{1}{2h}\begin{pmatrix}
            0 & 0 & 0 & 0 & 1
        \end{pmatrix}^T\otimes \left(I_2\right) \\
        B_q &= \frac{1}{2h}\begin{pmatrix}
            1 & 0 & 0 & 0 & 0
        \end{pmatrix}^T\otimes \left(-\mathcal{P}_1^T - h\mathcal{P}_0^T\right)
    \end{aligned}
\]

\begin{figure}
    \centering
    \includegraphics[width=.98\linewidth]{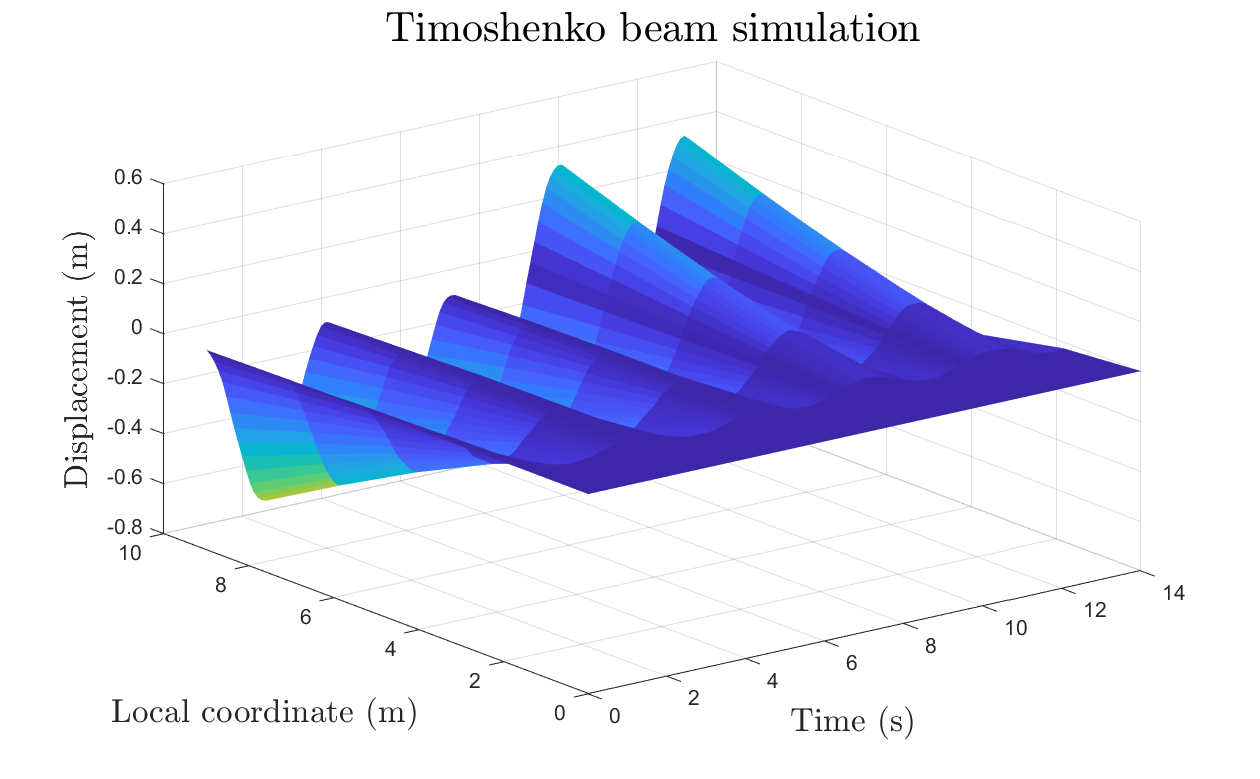}
    \caption{Timoshenko beam simulation.}
    \label{fig:TimoshenkoSimulation}
\end{figure}

Fig. \ref{fig:TimoshenkoSimulation} presents a simulation based on the finite-dimensional model, with time discretization performed using the midpoint method. The simulation shows the deformation of an aluminium rod (Young Modulus $E=68\,GPa$, Poisson's ratio $\nu=0.36$ and density $\rho=2698.9\,kg/m^3$), with a disc section of radius $2\,cm$ charged with a $2\,kg$ mass at the tip that is released at $7\,s$. In Fig \ref{fig:energyConservation} we can see that the energy is conserved when the inputs are set to 0.

\subsection{Example: Mindlin plate application}
For the 2D case, we present the discretization of the Mindlin plate equation given in PH form in \cite{ponce2024}. More precisely, we consider a rectangular plate clamped along one edge, with the two opposite edges free and the remaining edge supporting a mass suspended from a single off-center point. The interconnection matrices involve only neighboring points, as in \eqref{ec:I1_p_example}; and if the point is ahead in $\xi_i$, the connecting point is assigned a positive value; otherwise, it is negative. For $\mathcal{I}_0$, since there is no $\xi_0$ coordinate, all values are positive, as illustrated in \eqref{ec:I0_p_example}.

\begin{figure}[ht]
    \centering
    \includegraphics[trim={1cm 0.2cm 1cm 1cm},clip,width=.98\linewidth]{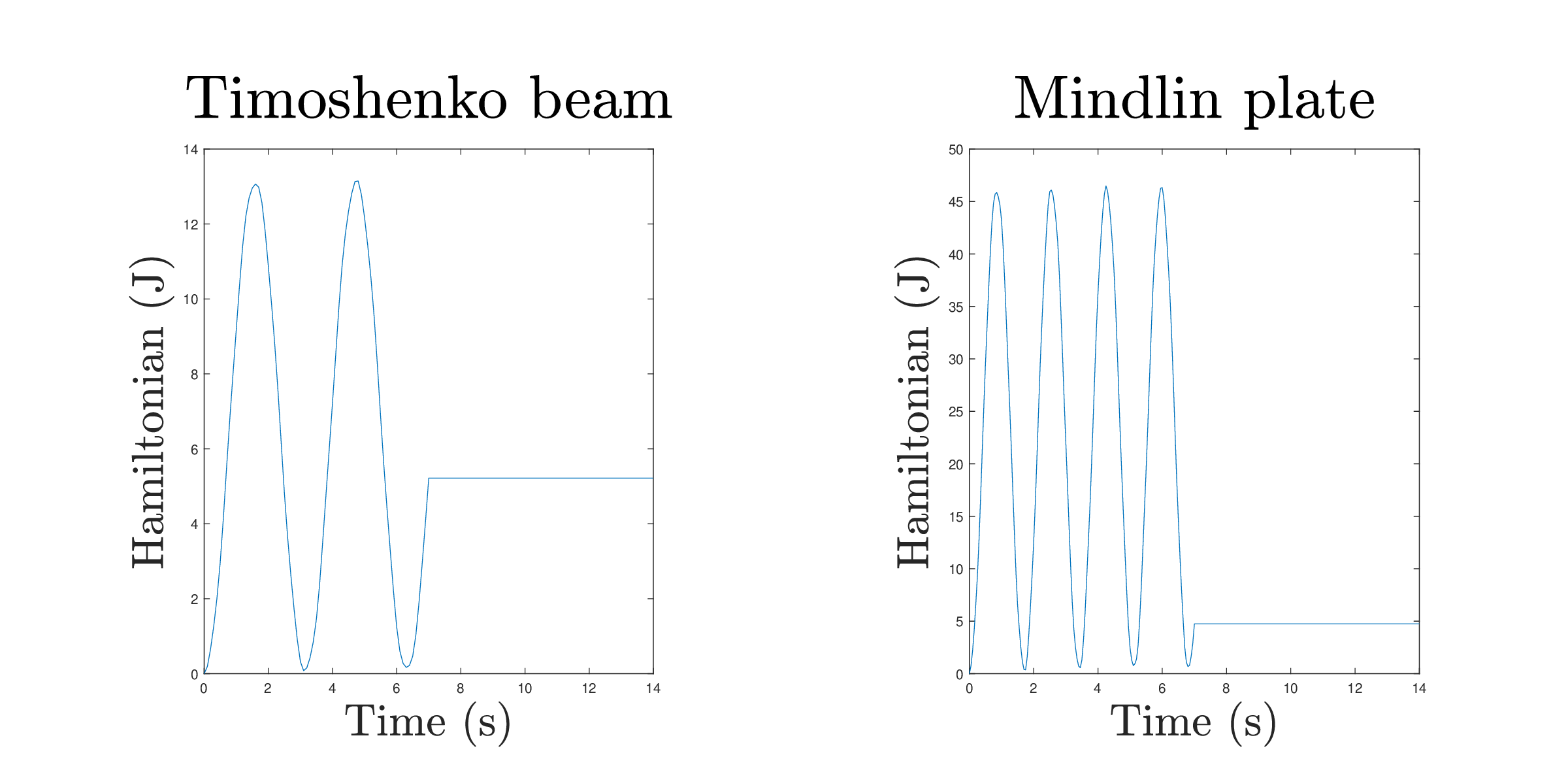}
    \caption{Hamiltonian over time for 1D and 2D simulations.}
    \label{fig:energyConservation}
\end{figure}

Fig. \ref{fig:energyConservation} shows the energy from a simulation conducted in the same manner as the previous one. The simulation considers an aluminum plate $3\,mm$ wide, supporting a $2\,kg$ mass that is released at $7\,s$. As in the 1D case, the energy is conserved when the mass is released.
\footnote{Due to the size limitation, animations of the simulations are included in \tt{https://github.com/I-DiazAl/ECC26-StaggeredGrid-PHS} }

\section{Conclusions}
In this paper we build on earlier results for structure-preserving discretization of 1D PH systems via staggered-grid finite differences, extending the approach to a more general class of 1D and 2D PH systems.  The proposed framework can be applied directly to 1D PH systems in which the Hamiltonian can be decomposed into two generalized coordinates and the interconnection operator is linear. Under similar conditions, the approach can also be applied to 2D systems, though only for specific spatial domains. The resulting finite-dimensional models preserve the PH structure and have explicit formulations, regardless of the input configuration.

Several extensions are currently being explored. One direction involves incorporating specific nonlinearities, either in the interconnection operator or via distributed inputs. Another is to adapt the discretization to triangular grids, which would enable the treatment of more general geometries. Finally, we are considering how to handle interconnections between systems that do not share exact port points. This would make the framework applicable to a wider range of spatial domains.

\bibliography{bibliography}

\end{document}